\documentclass[reqno]{amsart}

\newtheorem{Theorem}{Theorem}[section]
\newtheorem{Proposition}[Theorem]{Proposition}
\newtheorem{Corollary}[Theorem]{Corollary}
\newtheorem{Lemma}[Theorem]{Lemma}
\theoremstyle{remark}
\newtheorem{Remark}[Theorem]{Remark}
\numberwithin{equation}{section}

\allowdisplaybreaks

\begin{document}

\title[A nonterminating $_8\phi_7$ summation for $C_r$]
{A nonterminating
$\hbox{$_{\boldsymbol 8}\boldsymbol\phi_{\boldsymbol 7}$}$
summation for the root system $\hbox{$\boldsymbol C_{\boldsymbol r}$}$}
\author{Michael Schlosser}

\address{Institut f\"ur Mathematik der Universit\"at Wien,
Strudlhofgasse 4, A-1090 Wien, Austria}
\email{schlosse@ap.univie.ac.at}
\urladdr{http://www.mat.univie.ac.at/{\textasciitilde}schlosse}
\thanks{The author was supported by an APART grant of the Austrian
Academy of Sciences}
\date{November 8, 2002}


\begin{abstract}
Using multiple $q$-integrals and a determinant evaluation,
we establish a nonterminating $_8\phi_7$ summation for the
root system $C_r$. We also give some important specializations explicitly.
\end{abstract}

\maketitle

\section{Introduction}\label{sec0}
Bailey's~\cite[Eq.~(3.3)]{bailwp} nonterminating
very-well-poised $_8\phi_7$ summation,
\begin{multline}\label{87ntgl}
{}_8\phi_7\!\left[\begin{matrix}a,\,q\sqrt{a},-q\sqrt{a},b,c,d,e,f\\
\sqrt{a},-\sqrt{a},aq/b,aq/c,aq/d,aq/e,aq/f\end{matrix}\,;q,q\right]\\
+\frac{(aq,c,d,e,f,b/a,bq/c,bq/d,bq/e,bq/f;q)_{\infty}}
{(a/b,aq/c,aq/d,aq/e,aq/f,bc/a,bd/a,be/a,bf/a,b^2q/a;q)_{\infty}}\\\times
{}_8\phi_7\!\left[\begin{matrix}b^2/a,\,qb/\sqrt{a},-qb/\sqrt{a},
b,bc/a,bd/a,be/a,bf/a\\
b/\sqrt{a},-b/\sqrt{a},bq/a,bq/c,bq/d,bq/e,bq/f\end{matrix}\,;q,q\right]\\
=\frac{(aq,b/a,aq/cd,aq/ce,aq/cf,aq/de,aq/df,aq/ef;q)_{\infty}}
{(aq/c,aq/d,aq/e,aq/f,bc/a,bd/a,be/a,bf/a;q)_{\infty}},
\end{multline}
where $a^2q=bcdef$ (cf.~\cite[Eq.~(2.11.7)]{grhyp}), is one of the
deepest results in the classical theory of basic hypergeometric series.
It contains many important identities as special cases (such as
the nonterminating $_3\phi_2$ summation, the terminating $_8\phi_7$
summation, and all their specializations including the $q$-binomial theorem).
One way to derive \eqref{87ntgl} is to start with a particular rational
function identity, namely Bailey's~\cite{bail10} very-well-poised
$_{10}\phi_9$ transformation, and apply a nontrivial limit procedure,
see the exposition in Gasper and Rahman~\cite[Secs.~2.10 and 2.11]{grhyp}.

Basic hypergeometric series (and, more generally, $q$-series)
have various applications in combinatorics, number theory,
representation theory, statistics, and physics,
see Andrews~\cite{andappl}, \cite{qandrews}.  For a general account
of the importance of basic hypergeometric series in the theory of
special functions see Andrews, Askey, and Roy~\cite{sfaar}.

There are different types of multivariable series. The one we are
concerned with are so-called {\em multiple basic hypergeometric series
associated to root systems} (or, equivalently, to {\em Lie algebras}).
This is mainly just a classification of certain multiple series according
to the type of specific factors (such as a Vandermonde determinant)
appearing in the summand. We omit giving a precise definition here,
but instead refer to papers of Bhatnagar~\cite{bhatdn} or
Milne~\cite[Sec.~5]{milnetf}.

The significance of the nonterminating $_8\phi_7$ summation
\eqref{87ntgl} lies in the fact that it can be used for deriving
other nonterminating transformation formulae, see Gasper and
Rahman~\cite[Secs.~ 2.12 and 3.8]{grhyp}, and Schlosser~\cite{schleltf}.
Thus, it is apparently desirable to find (various) multivariable
generalizations of Bailey's nonterminating $_8\phi_7$ summation.

In this paper, we give a multivariable nonterminating $_8\phi_7$
summation for the root system $C_r$ (or, equivalently,
the symplectic group $Sp(r)$), see Corollary~\ref{cnnt87}.
We deduce this result from an equivalent multiple $q$-integral
evaluation, Theorem~\ref{cnint87}. 
In our proof of the latter we utilize a simple determinant method,
essentially the same which was introduced by Gustafson and
Krattenthaler~\cite{guskradet} and which we further exploited
in \cite{schlhypdet} to derive a number of identities
for multiple basic hypergeometric series.
The difference here is that now we apply the method to
{\em integrals} and {\em $q$-integrals} whereas in
\cite{schlhypdet} we had only applied it to sums.
Our new $C_r$ nonterminating $_8\phi_7$ summation is not the first
multivariable nonterminating $_8\phi_7$ sum that has been found.
In fact, Degenhardt and Milne~\cite{degmilne} already derived
such a result for the root system $A_n$ (or, equivalently,
the unitary group $U\!(n)$), a result we consider to be
deeper than ours.
While Corollary~\ref{cnnt87} is derived by elementary means,
by combining known one-variable results with
the argument of interchanging the order of summations,
or of summation and ($q$-)integration (this is what the
determinant method in this article really does),
Degenhardt and Milne deduce their multivariable summation formula
by extending Gasper and Rahman's~\cite[Sec.~ 2.10]{grhyp} analysis
to higher dimensions which appears to be fairly nontrivial.
However, supported by the combinatorial applications
(in Krattenthaler~\cite{krattmaj}, and Gessel and
Krattenthaler~\cite{geskratt}) of identities of type
strikingly similar to the one being investigated in this paper,
we believe that the identities derived here will have
future applications and deserve being written out in detail.

Our paper is organized as follows: In Section~\ref{secpre},
we review some basics in the theory of basic hypergeometric series.
Further, we also note a determinant lemma which we need as an
ingredient in proving our results in Sections~\ref{secmbeta} and
\ref{secmain}. We demonstrate the method of proof in
Section~\ref{secmbeta} by deriving a simple multidimensional beta
integral evaluation. In Section~\ref{secmain}, we derive an
(attractive) multiple $q$-integral evauation, Theorem~\ref{cnint87},
which in Section~\ref{secm87} is used to explicitly write out a
nonterminating $_8\phi_7$ summation for the root system $C_r$,
see Corollary~\ref{cnnt87}.
Finally, in Section~\ref{secspec} we explicitly list several
interesting specializations of Theorem~\ref{cnint87}
(and of the equivalent Corollary~\ref{cnnt87}).

\section{Basic hypergeometric series and a determinant lemma}\label{secpre}

Here we recall some standard notation for $q$-series,
and basic hypergeometric series (cf.\ \cite{grhyp}).

Let $q$ be a complex number such that $0<|q|<1$. We define the
{\em $q$-shifted factorial} for all integers $k$ by
\begin{equation*}
(a;q)_{\infty}:=\prod_{j=0}^{\infty}(1-aq^j)\qquad\text{and}\qquad
(a;q)_k:=\frac{(a;q)_{\infty}}{(aq^k;q)_{\infty}}.
\end{equation*}
For brevity, we employ the condensed notation
\begin{equation*}
(a_1,\ldots,a_m;q)_k\equiv (a_1;q)_k\dots(a_m;q)_k
\end{equation*}
where $k$ is an integer or infinity. Further, we utilize
\begin{equation}\label{defhyp}
{}_s\phi_{s-1}\!\left[\begin{matrix}a_1,a_2,\dots,a_s\\
b_1,b_2,\dots,b_{s-1}\end{matrix}\,;q,z\right]:=
\sum _{k=0} ^{\infty}\frac {(a_1,a_2,\dots,a_s;q)_k}
{(q,b_1,\dots,b_{s-1};q)_k}z^k,
\end{equation}
to denote the {\em basic hypergeometric ${}_s\phi_{s-1}$ series}.
In \eqref{defhyp}, $a_1,\dots,a_s$ are
called the {\em upper parameters}, $b_1,\dots,b_{s-1}$ the
{\em lower parameters}, $z$ is the {\em argument}, and 
$q$ the {\em base} of the series.
The series in \eqref{defhyp} terminates if one of the upper parameters,
say $a_s$, equals $q^{-n}$ where $n$ is a nonnegeative integer.
If the series does not terminate, we need $|z|<1$ for convergence. 

The classical theory of basic hypergeometric series contains
numerous summation and transformation formulae
involving ${}_s\phi_{s-1}$ series.
Many of these summation theorems require
that the parameters satisfy the condition of being
either balanced and/or very-well-poised.
An ${}_s\phi_{s-1}$ basic hypergeometric series is called
{\em balanced} if $b_1\cdots b_{s-1}=a_1\cdots a_sq$ and $z=q$.
An ${}_s\phi_{s-1}$ series is {\em well-poised} if
$a_1q=a_2b_1=\cdots=a_sb_{s-1}$ and {\em very-well-poised}
if it is well-poised and if $a_2=-a_3=q\sqrt{a_1}$.
Note that the factor
\begin{equation*}
\frac {1-a_1q^{2k}}{1-a_1}
\end{equation*}
appears in  a very-well-poised series.
The parameter $a_1$ is usually referred to as the
{\em special parameter} of such a series.

One of the most important summation theorems in the theory of basic
hypergeometric series
is Jackson's~\cite{jacksum} terminating very-well-poised balanced
${}_8\phi_7$ summation (cf.\ \cite[Eq.~(2.6.2)]{grhyp}):
\begin{multline}\label{87gl}
{}_8\phi_7\!\left[\begin{matrix}a,\,q\sqrt{a},-q\sqrt{a},b,c,d,
a^2q^{1+n}/bcd,q^{-n}\\
\sqrt{a},-\sqrt{a},aq/b,aq/c,aq/d,bcdq^{-n}/a,aq^{1+n}\end{matrix}\,;q,
q\right]\\
=\frac {(aq,aq/bc,aq/bd,aq/cd;q)_n}
{(aq/b,aq/c,aq/d,aq/bcd;q)_n}.
\end{multline}
Clearly, \eqref{87gl} is the special case $f\to q^{-n}$ of \eqref{87ntgl}.

For studying nonterminating basic hypergeometric series it is often
convenient to utilize Jackson's~\cite{jackson} $q$-integral notation,
defined by
\begin{equation}\label{qint}
\int_a^bf(t)\,d_qt=\int_0^bf(t)\;d_qt-\int_0^af(t)\;d_qt,
\end{equation}
where
\begin{equation}\label{qint0a}
\int_0^af(t)\;d_qt=a(1-q)\sum_{k=0}^{\infty}f(aq^k)q^k.
\end{equation}
If $f$ is continuous on $[0,a]$, then it is easily seen that
\begin{equation*}
\lim_{q\to1^-}\int_0^a f(t)\;d_qt=\int_0^a f(t)\;dt,
\end{equation*}
see~\cite[Eq.~(1.11.6)]{grhyp}.

Using the above $q$-integral notation, the nonterminating $_8\phi_7$
summation \eqref{87ntgl} can be conveniently expressed as
\begin{multline}\label{int87gl}
\int_a^b\frac{(qt/a,qt/b,t/\sqrt{a},-t/\sqrt{a},qt/c,qt/d,qt/e,qt/f;
q)_{\infty}}{(t,bt/a,qt/\sqrt{a},-qt/\sqrt{a},ct/a,dt/a,et/a,ft/a;
q)_{\infty}}\;d_qt\\
=\frac{b(1-q)(q,a/b,bq/a,aq/cd,aq/ce,aq/cf,aq/de,aq/df,aq/ef;q)_{\infty}}
{(b,c,d,e,f,bc/a,bd/a,be/a,bf/a;q)_{\infty}},
\end{multline}
where $a^2q=bcdef$ (cf.\ \cite[Eq.~(2.11.8)]{grhyp}).

A standard reference for basic hypergeometric series
is Gasper and Rahman's text~\cite{grhyp}.
In our computations in the subsequent sections
we frequently use some elementary identities of
$q$-shifted factorials, listed in \cite[Appendix~I]{grhyp}.

The following determinant evaluation was given as Lemma~A.1 in
\cite{schlhypdet} where it was derived from a determinant lemma
of Krattenthaler~\cite[Lemma~34]{krattmaj}.

\begin{Lemma}\label{lemdet1}
Let $X_1,\dots,X_r$, $A$, $B$, and $C$ be indeterminate. Then there holds
\begin{multline}\label{lemdet1gl}
\det_{1\le i,j\le r}\left(
\frac{(AX_i,AC/X_i;q)_{r-j}}
{(BX_i,BC/X_i;q)_{r-j}}\right)=
\prod_{1\le i<j\le r}(X_j-X_i)(1-C/X_iX_j)\\\times
A^{\binom r2}q^{\binom r3}
\prod_{i=1}^{r}\frac {(B/A,ABCq^{2r-2i};q)_{i-1}}
{(BX_i,BC/X_i;q)_{r-1}}.
\end{multline}
\end{Lemma}

The above determinant evaluation was generalized to the elliptic case
(more precisely, to an evaluation involving Jacobi theta functions)
by Warnaar~\cite[Cor.~5.4]{warnell}.

\section{A multidimensional beta integral evaluation}\label{secmbeta}

Here, we present a simple multivariable extension of Euler's
beta inetgral evaluation. The proof serves as an illustration
of the determinant method which we use in Section \ref{secmain}
to derive a multivariable extension of \eqref{int87gl}.

\begin{Proposition}\label{mbetaint}
Let $a$, $b$, and $x_1,\dots,x_r$ be indeterminate. Then
\begin{multline}\label{mbetaintgl}
\int_0^1\dots\int_0^1\prod_{1\le i<j\le r}(u_i-u_j)
\prod_{i=1}^ru_i^{a-1+x_i}(1-u_i)^{b-1}\,du_r\dots du_1\\
=\prod_{1\le i<j\le r}(x_i-x_j)\prod_{i=1}^r
\frac{\Gamma(a+x_i)\,\Gamma(b+i-1)}{\Gamma(a+b+x_i+r-1)},
\end{multline}
provided $\Re (a+x_i), \Re(b)>0$, for $i=1,\dots,r$.
\end{Proposition}

\begin{Remark}
We note the differences between Proposition~\ref{mbetaint} 
and Selberg's~\cite{selberg} integral,
\begin{multline}\label{selgl}
\int_0^1\dots\int_0^1\prod_{1\le i<j\le r}|u_i-u_j|^{2c}
\prod_{i=1}^ru_i^{a-1}(1-u_i)^{b-1}\,du_r\dots du_1\\
=\prod_{i=1}^r
\frac{\Gamma(a+(i-1)c)\,\Gamma(b+(i-1)c)\,\Gamma(ic+1)}
{\Gamma(a+b+(r+i-2)c)\,\Gamma(c+1)},
\end{multline}
where $\Re (a),\Re (b)>0$, and
$\Re (c)>\max(-1/r,-\Re (a)/(r-1),-\Re (b)/(r-1))$. In \eqref{mbetaintgl}
we have additional parameters $x_1,\dots,x_r$, while in
\eqref{selgl} the absolute value of the discriminant
$\prod_{1\le i<j\le r}(u_i-u_j)$ in the integrand is taken to an
arbritrary power $2c$, which makes the computation considerably
more difficult. 
\end{Remark}

\begin{proof}[Proof of Proposition~\ref{mbetaint}]
First, we note that if in Lemma~\ref{lemdet1}
we let $C\to 0$, replace $A$, $B$, and $X_i$ by
$q^a$, $q^{a+b}$, and $q^{x_i}$, for $i=1,\dots,r$, respectively,
and then let $q\to 1$, we have
\begin{equation}\label{detcor}
\det_{1\le i,j\le r}\left(\frac{(a+x_i)_{r-j}}{(a+b+x_i)_{r-j}}\right)
=\prod_{1\le i<j\le r}(x_i-x_j)
\prod_{i=1}^{r}\frac {(b)_{i-1}}
{(a+b+x_i)_{r-1}},
\end{equation}
where
\begin{equation}\label{poch}
(\alpha)_k:=\frac{\Gamma(a+k)}{\Gamma(a)}
\end{equation}
is the {\em shifted factorial}.
Thus,
\begin{multline*}
\int_0^1\dots\int_0^1\prod_{1\le i<j\le r}(u_i-u_j)
\prod_{i=1}^ru_i^{a-1+x_i}(1-u_i)^{b-1}\,du_r\dots du_1\\
=\det_{1\le i,j\le r}\left(\int_0^1u_i^{a-1+x_i+r-j}
(1-u_i)^{b-1}\,du_i\right)\\
=\det_{1\le i,j\le r}\left(\frac
{\Gamma(a+x_i+r-j)\,\Gamma(b)}{\Gamma(a+b+x_i+r-j)}\right)
=\det_{1\le i,j\le r}\left(\frac
{\Gamma(a+x_i)\,\Gamma(b)}{\Gamma(a+b+x_i)}
\frac{(a+x_i)_{r-j}}{(a+b+x_i)_{r-j}}\right)\\
=\prod_{i=1}^r\frac
{\Gamma(a+x_i)\,\Gamma(b)}{\Gamma(a+b+x_i)}\cdot
\det_{1\le i,j\le r}\left(
\frac{(a+x_i)_{r-j}}{(a+b+x_i)_{r-j}}\right)\\
=\prod_{i=1}^r\frac
{\Gamma(a+x_i)\,\Gamma(b)}{\Gamma(a+b+x_i)}
\prod_{1\le i<j\le r}(x_i-x_j)
\prod_{i=1}^{r}\frac {(b)_{i-1}}
{(a+b+x_i)_{r-1}}\\
=\prod_{1\le i<j\le r}(x_i-x_j)\prod_{i=1}^r
\frac{\Gamma(a+x_i)\,\Gamma(b+i-1)}{\Gamma(a+b+x_i+r-1)},
\end{multline*}
where we have used linearity of the determinant with respect to rows,
the Vandermonde determinant evaluation
$\det_{1\le i,j\le r}(u_i^{r-j})=\prod_{1\le i<j\le r}(u_i-u_j)$,
Euler's beta integral evaluation
$\int_0^1u^{a-1}(1-u)^{b-1}dt=\frac{\Gamma(a)\,\Gamma(b)}{\Gamma(a+b)}$,
for $\Re(a),\Re(b)>0$, the definition of the shifted factorial
\eqref{poch}, and the determinant evaluation \eqref{detcor}.
\end{proof}

\section{A multiple $q$-integral evaluation}\label{secmain}

By iteration, the extension of \eqref{qint0a} to {\em multiple}
$q$-integrals is straightforward:
\begin{multline}\label{mqint0a}
\int_0^{a_1}\dots\int_0^{a_r}f(t_1,\dots,t_r)\;d_qt_r\dots d_qt_1\\
=a_1\dots a_r(1-q)^r\sum_{k_1,\dots,k_r=0}^{\infty}f(a_1q^{k_1},\dots,
a_rq^{k_r})q^{k_1+\dots+k_r}.
\end{multline}
Similarly, the extension of \eqref{qint} is
\begin{multline}\label{mqint}
\int_{a_1}^{b_1}\dots\int_{a_r}^{b_r}f(t_1,\dots,t_r)\;d_qt_r\dots d_qt_1\\
=\sum_{S\subseteq\{1,2,\dots,r\}}\left(\prod_{i\in S}(-a_i)\right)
\left(\prod_{i\notin S}b_i\right)(1-q)^r\\\times
\sum_{k_1,\dots,k_r=0}^{\infty}f(c_1(S)q^{k_1},\dots,c_r(S)q^{k_r})
q^{k_1+\dots+k_r},
\end{multline}
where the outer sum runs over all $2^r$ subsets $S$ of $\{1,2,\dots,r\}$,
and where $c_i(S)=a_i$ if $i\in S$ and $c_i(S)=b_i$ if $i\notin S$,
for $i=1,\dots,r$.

We give our main result, a $C_r$ extension of \eqref{int87gl}:

\begin{Theorem}\label{cnint87}
Let $a^2q^{2-r}=bcdef$. Then there holds
\begin{multline}\label{cnint87gl}
\int_{ax_1}^{b}\dots\int_{ax_r}^{b}
\prod_{1\le i<j\le r}(t_i-t_j)(1-t_it_j/a)
\prod_{i=1}^r(1-t_i^2/a)\\
\times\prod_{i=1}^r
\frac{(qt_i/ax_i,qt_i/b,qt_i/c,qt_i/d,qt_i/e,qt_ix_i/f;q)_{\infty}}
{(t_ix_i,bt_i/a,ct_i/a,dt_i/a,et_i/a,ft_i/ax_i;q)_{\infty}}\;
d_qt_r\dots d_qt_1\\
=a^{\binom r2}b^r(1-q)^r
\prod_{1\le i<j\le r}(x_i-x_j)(1-ax_ix_j/f)\\\times
\prod_{i=1}^r\frac{(q,ax_i/b,bq/ax_i,aq^{2-i}/cd,aq^{2-i}/ce;q)_{\infty}}
{(bx_i,cx_i,dx_i,ex_i,f;q)_{\infty}}\\\times
\prod_{i=1}^r\frac{(ax_iq/cf,aq^{2-i}/de,ax_iq/df,ax_iq/ef;q)_{\infty}}
{(bcq^{i-1}/a,bdq^{i-1}/a,beq^{i-1}/a,bf/ax_i;q)_{\infty}}.
\end{multline}
\end{Theorem}

\begin{proof}
We have
\begin{multline*}
\prod_{1\le i<j\le r}(t_i-t_j)(1-t_it_j/a)
=\prod_{i=1}^r\frac{(q^{2-r}t_i/d,aq^{2-r}/dt_i;q)_{r-1}}
{(aq^{2-r}/cd,cq^{2+r-2i}/d;q)_{i-1}}t_i^{r-1}\\\times
c^{-\binom r2}q^{-\binom r3}\;
\det_{1\le i,j\le r}\left(\frac {(ct_i/a,c/t_i;q)_{r-j}}
{(q^{2-r}t_i/d,aq^{2-r}/dt_i;q)_{r-j}}\right),
\end{multline*}
due to the $X_i\mapsto t_i$, $A\mapsto c/a$, $B\mapsto q^{2-r}/d$,
and $C\mapsto a$ case of Lemma~\ref{lemdet1}.
Hence, using some elementary identities from \cite[Appendix~I]{grhyp},
we may write the left hand side of \eqref{cnint87gl} as
\begin{multline*}
\left(\frac ad\right)^{\binom r2}q^{-\binom r3}
\prod_{i=1}^r(aq^{2-r}/cd,cq^{2+r-2i}/d;q)_{i-1}^{-1}\\\times
\det_{1\le i,j\le r}\Bigg(\int_{ax_i}^b
\frac {(1-t_i^2/a)(t_iq/ax_i,t_iq/b;q)_{\infty}}
{(t_ix_i,bt_i/a;q)_{\infty}}\\\times
\frac {(t_iq^{1-r+j}/c,t_iq^{2-j}/d,t_iq/e,t_ix_iq/f;q)_{\infty}}
{(ct_iq^{r-j}/a,dt_iq^{j-1}/a,et_i/a,ft_i/ax_i;q)_{\infty}}\;
d_qt_i\Bigg).
\end{multline*}
Now, to the integral inside the determinant we apply the
$q$-integral evaluation \eqref{int87gl},
with the substitution $t\mapsto t_ix_i$, and the
replacements $a\mapsto ax_i^2$, $b\mapsto bx_i$, $c\mapsto cq^{r-j}x_i$,
$d\mapsto dq^{j-1}x_i$, and $e\mapsto ex_i$. Thus we obtain
\begin{multline*}
\left(\frac ad\right)^{\binom r2}q^{-\binom r3}
\prod_{i=1}^r(aq^{2-r}/cd,cq^{2+r-2i}/d;q)_{i-1}^{-1}\\\times
\det_{1\le i,j\le r}\Bigg(
\frac{b(1-q)(q,ax_i/b,bq/ax_i,aq^{2-r}/cd,aq^{1-r+j}/ce;q)_{\infty}}
{(bx_i,cx_iq^{r-j},dx_iq^{j-1},ex_i,f;q)_{\infty}}\\\times
\frac{(ax_iq^{1-r+j}/cf,aq^{2-j}/de,ax_iq^{2-j}/df,ax_iq/ef;q)_{\infty}}
{(bcq^{r-j}/a,bdq^{j-1}/a,be/a,bf/ax_i;q)_{\infty}}
\Bigg).
\end{multline*}
Now, by using linearity of the determinant with respect to rows and columns,
we take some factors out of the determinant and obtain
\begin{multline*}
\left(\frac ad\right)^{\binom r2}q^{-\binom r3}b^r(1-q)^r
\prod_{i=1}^r\frac{(q,ax_i/b,bq/ax_i;q)_{\infty}}
{(aq^{2-r}/cd,cq^{2+r-2i}/d;q)_{i-1}}\\\times
\prod_{i=1}^r
\frac{(aq^{2-r}/cd,aq^{1-r+i}/ce,ax_iq/cf,aq^{2-i}/de,
ax_iq^{2-r}/df,ax_iq/ef;q)_{\infty}}
{(bx_i,cx_i,dx_iq^{r-1},ex_i,f,bcq^{r-i}/a,bdq^{i-1}/a,
be/a,bf/ax_i;q)_{\infty}}\\\times
\left(\frac a{cdf}\right)^{\binom r2}q^{-3\binom r3}
\det_{1\le i,j\le r}\Bigg(
\frac{(cx_i,cf/ax_i;q)_{r-j}}
{(ax_iq^{2-r}/df,q^{2-r}/dx_i;q)_{r-j}}
\Bigg).
\end{multline*}
The determinant can be evaluated by means of Lemma~\ref{lemdet1} with
$X_i\mapsto x_i$, $A\mapsto c$, $B\mapsto aq^{2-r}/df$, and
$C\mapsto f/a$; specifically
\begin{multline*}
\det_{1\le i,j\le r}\Bigg(
\frac{(cx_i,cf/ax_i;q)_{r-j}}
{(ax_iq^{2-r}/df,q^{2-r}/dx_i;q)_{r-j}}\Bigg)\\=
c^{\binom r2}q^{\binom r3}
\prod_{1\le i<j\le r}(x_j-x_i)(1-f/ax_ix_j)
\prod_{i=1}^r\frac{(aq^{2-r}/cdf,cq^{2+r-2i}/d;q)_{i-1}}
{(ax_iq^{2-r}/df,q^{2-r}/dx_i;q)_{r-1}}.
\end{multline*}
Substituting our calculations and performing further elementary
manipulations we arrive at the right hand side of \eqref{cnint87gl}. 
\end{proof}

\section{A multivariable nonterminating ${}_8\phi_7$ summation}\label{secm87}

Note that if the integrand $f(t_1,\ldots,t_r)$ of the multiple integral
in \eqref{cnint87gl} were an antisymmetric function in $t_1,\dots,t_r$,
the multiple sum in \eqref{mqint} would simplify considerably.
In fact, if $t_i=bq^{k_i}$ and $t_j=bq^{k_j}$, for a pair $i<j$,
we would then have
\begin{equation*}
\sum_{k_i,k_j=0}^{\infty}f(\ldots,bq^{k_i},\ldots,bq^{k_j},\ldots)=0.
\end{equation*}
(A double sum of any function antisymmetric in its two summation indices
vanishes.)
Thus, the multiple $q$-integral
$\int_{ax_1}^b\dots\int_{ax_r}^b f(t_1,\ldots,t_r)\;d_qt_r\dots d_qt_1$,
being a sum of $2^r$ sums according to \eqref{mqint}, would reduce
to a sum of $r+1$ nonzero sums. In particular, we would have
\begin{multline}
\int_{ax_1}^b\dots\int_{ax_r}^b f(t_1,\ldots,t_r)\;d_qt_r\dots d_qt_1\\
=(-1)^ra^rx_1\dots x_r(1-q)^r
\sum_{k_1,\dots,k_r=0}^{\infty}f(ax_1q^{k_1},\dots,ax_rq^{k_r})\,
q^{\sum_{i=1}^r k_i}\\
+(-1)^{r-1}a^{r-1}bx_1\dots x_r(1-q)^r\sum_{l=1}^rx_l^{-1}\\\times
\sum_{k_1,\dots,k_r=0}^{\infty}
f(ax_1q^{k_1},\dots,ax_{l-1}q^{k_{l-1}},bq^{k_l},
ax_{l+1}q^{k_{l+1}},\dots,ax_rq^{k_r})\,q^{\sum_{i=1}^r k_i}.
\end{multline}
A very similar situation occurs in the $U\!(n)$ (or $A_r$) nonterminating
$_8\phi_7$ summation by Degenhardt and Milne~\cite{degmilne}
(however, their argument is reversed, i.e., they first derive a
nonterminating summation and then deduce the multiple $q$-integral
evaluation).
Unfortunately, in our case the multiple integrand in \eqref{cnint87gl}
is {\em not} antisymmetric in $t_1,\dots,t_r$, whence we have all
$2^r$ sums on the right hand side of \eqref{mqint}.

We write out \eqref{mqint} explicitly for our integral in \eqref{cnint87gl}:
\begin{multline*}
\int_{ax_1}^{b}\dots\int_{ax_r}^{b}
\prod_{1\le i<j\le r}(t_i-t_j)(1-t_it_j/a)
\prod_{i=1}^r(1-t_i^2/a)\\
\times\prod_{i=1}^r
\frac{(qt_i/ax_i,qt_i/b,qt_i/c,qt_i/d,qt_i/e,qt_ix_i/f;q)_{\infty}}
{(t_ix_i,bt_i/a,ct_i/a,dt_i/a,et_i/a,ft_i/ax_i;q)_{\infty}}\;
d_qt_r\dots d_qt_1\\
=\sum_{S\subseteq\{1,2,\dots,r\}}
(-1)^{|S|}a^{|S|}b^{r-|S|}(1-q)^ra^{\binom{|S|}2}b^{\binom{r-|S|}2}
\prod_{i\in S}x_i\\\times
\sum_{k_1,\dots,k_r=0}^{\infty}
\underset{i,j\in S}{\prod_{1\le i<j\le r}}
(x_iq^{k_i}-x_jq^{k_j})(1-ax_ix_jq^{k_i+k_j})
\prod_{i\in S}(1-ax_i^2q^{2k_i})\\\times
\underset{i,j\notin S}{\prod_{1\le i<j\le r}}
(q^{k_i}-q^{k_j})(1-b^2q^{k_i+k_j}/a)
\prod_{i\notin S}(1-bq^{2k_i}/a)\\\times
\prod_{i\in S, j\notin S}(ax_iq^{k_i}-bq^{k_j})(1-bx_ix_jq^{k_i+k_j})
(-1)^{\chi(i>j)}\\\times
\prod_{i\in S}
\frac{(q^{1+k_i},ax_iq^{1+k_i}/b,ax_iq^{1+k_i}/c,ax_iq^{1+k_i}/d,
ax_iq^{1+k_i}/e,ax_i^2q^{1+k_i}/f;q)_{\infty}}
{(ax_i^2q^{k_i},bx_iq^{k_i},cx_iq^{k_i},dx_iq^{k_i},ex_iq^{k_i},
fq^{k_i};q)_{\infty}}\\\times
\prod_{i\notin S}
\frac{(bq^{1+k_i}/ax_i,q^{1+k_i},bq^{1+k_i}/c,bq^{1+k_i}/d,
bq^{1+k_i}/e,bx_iq^{1+k_i}/f;q)_{\infty}}
{(bx_iq^{k_i},b^2q^{k_i}/a,bcq^{k_i}/a,bdq^{k_i}/a,beq^{k_i}/a,
bfq^{k_i}/ax_i;q)_{\infty}}\;\cdot q^{\sum_{i=1}^r k_i},
\end{multline*}
where $|S|$ denotes the number of elements of $S$, and $\chi$
is the truth function (which evaluates to one if the argument is true
and evaluates to zero otherwise).
Now, if we set the obtained sum of $2^r$ sums equal to the right hand side of
\eqref{cnint87gl} and divide both sides by
\begin{multline*}
(-1)^ra^{\binom{r+1}2}x_1\dots x_r(1-q)^r
\prod_{1\le i<j\le r}(x_i-x_j)(1-ax_ix_j)\\\times
\prod_{i=1}^r\frac{(q,ax_iq/b,ax_iq/c,ax_iq/d,ax_iq/e,ax_i^2q/f;q)_{\infty}}
{(ax_i^2q,bx_i,cx_i,dx_i,ex_i,f;q)_{\infty}},
\end{multline*}
and simplify, we obtain the following result which reduces to
\eqref{87ntgl} when $r=1$.

\begin{Corollary}[A $C_r$ nonterminating ${}_8\phi_7$ summation]\label{cnnt87}
Let $a^2q^{2-r}=bcdef$. Then there holds
\begin{multline}\label{cnnt87gl}
\sum_{S\subseteq\{1,2,\dots,r\}}
\left(\frac ba\right)^{\binom{r-|S|}2}\prod_{i\notin S}
\frac{(ax_i^2q,cx_i,dx_i,ex_i,f;q)_{\infty}}
{(ax_i/b,ax_iq/c,ax_iq/d,ax_iq/e,ax_i^2q/f;q)_{\infty}}\\\times
\prod_{i\notin S}
\frac{(b/ax_i,bq/c,bq/d,bq/e,bx_iq/f;q)_{\infty}}
{(b^2q/a,bc/a,bd/a,be/a,bf/ax_i;q)_{\infty}}\\\times
\sum_{k_1,\dots,k_r=0}^{\infty}
\underset{i,j\in S}{\prod_{1\le i<j\le r}}
\frac{(x_iq^{k_i}-x_jq^{k_j})(1-ax_ix_jq^{k_i+k_j})}
{(x_i-x_j)(1-ax_ix_j)}
\prod_{i\in S}\frac{(1-ax_i^2q^{2k_i})}{(1-ax_i^2)}\\\times
\underset{i,j\notin S}{\prod_{1\le i<j\le r}}
\frac{(q^{k_i}-q^{k_j})(1-b^2q^{k_i+k_j}/a)}
{(x_i-x_j)(1-ax_ix_j)}
\prod_{i\notin S}\frac{(1-bq^{2k_i}/a)}{(1-b^2/a)}\\\times
\prod_{i\in S, j\notin S}
\frac{(x_iq^{k_i}-bq^{k_j}/a)(1-bx_ix_jq^{k_i+k_j})}
{(x_i-x_j)(1-ax_ix_j)}\\\times
\prod_{i\in S}\frac{(ax_i^2,bx_i,cx_i,dx_i,ex_i,f;q)_{k_i}}
{(q,ax_iq/b,ax_iq/c,ax_iq/d,ax_iq/e,ax_i^2q/f;q)_{k_i}}\\\times
\prod_{i\notin S}\frac{(b^2/a,bx_i,bc/a,bd/a,be/a,bf/ax_i;q)_{k_i}}
{(q,bq/ax_i,bq/c,bq/d,bq/e,bx_iq/f;q)_{k_i}}\;\cdot q^{\sum_{i=1}^r k_i}\\
=\prod_{1\le i<j\le r}\frac{(1-ax_ix_j/f)}{(1-ax_ix_j)}
\prod_{i=1}^r\frac{(ax_i^2q,b/ax_i,aq^{2-i}/cd,aq^{2-i}/ce;q)_{\infty}}
{(ax_iq/c,ax_iq/d,ax_iq/e,ax_i^2q/f;q)_{\infty}}\\\times
\prod_{i=1}^r\frac{(ax_iq/cf,aq^{2-i}/de,ax_iq/df,ax_iq/ef;q)_{\infty}}
{(bcq^{i-1}/a,bdq^{i-1}/a,beq^{i-1}/a,bf/ax_i;q)_{\infty}}.
\end{multline}
\end{Corollary}

\section{Specializations}\label{secspec}

It is clear that in Corollary~\ref{cnnt87}, if we replace
$e$ by $a^2q^{2-r}/bcdf$ and then let $f\to q^{-N}$, we obtain
a $C_r$ extension of Jackson's $_8\phi_7$ summation \eqref{87gl}.
This result,
\begin{multline}\label{phi87gl}
\sum_{k_1,\dots,k_r=0}^N
\prod_{1\le i<j\le r}\frac{(x_iq^{k_i}-x_jq^{k_j})
(1-a x_i x_j q^{k_i+k_j})}{(x_i-x_j)(1-ax_ix_j)}
\prod_{i=1}^r\frac {(1-ax_i^2q^{2k_i})}{(1-ax_i^2)}\\
\times\prod_{i=1}^r\frac {(ax_i^2,bx_i,cx_i,dx_i,
a^2x_iq^{2-r+N}/bcd,q^{-N};q)_{k_i}}
{(q,ax_iq/b,ax_iq/c,ax_iq/d,bcdx_iq^{r-1-N}/a,
ax_i^2q^{1+N};q)_{k_i}}\,\cdot q^{\sum_{i=1} ^rk_i}\\
=\prod _{1\le i<j\le r}
\frac{(1-ax_ix_jq^N)}{(1-ax_ix_j)}
\prod_{i=1}^r\frac{(ax_i^2q,aq^{2-i}/bc,aq^{2-i}/bd,
aq^{2-i}/cd;q)_N}
{(aq^{2-r}/bcdx_i,ax_iq/d,ax_iq/c,ax_iq/b;q)_N},
\end{multline}
was given as Theorem~4.3 in \cite{schlhypdet}. An extension of
\eqref{phi87gl} to {\em elliptic} hypergeometric series
was found by Warnaar~\cite[Theorem~5.1]{warnell}.

Next, if in \eqref{cnint87gl},
we first replace $e$ by $a^2q^{2-r}/bcdf$, then 
$b$ by $aq/b$, do the substitution $t_i\mapsto at_i$,
divide both sides by $a^{\binom{r+1}2}$, let $a\to 0$,
and afterwards do the simultaneous substitutions
$x_i\mapsto x_i\sqrt{a}$, $b\mapsto q\sqrt{a}/b$,
$c\mapsto cq^{1-r}/bdf\sqrt{a}$, $d\mapsto d\sqrt{a}$,
$f\mapsto fa$, then multiply both sides by $\sqrt{a}^{\binom{r+1}2}$,
and perform the substitution $t_i\mapsto t_i/\sqrt{a}$, 
we obtain the following multiple $q$-integral evaluation which is an
$r$-dimensional extension of Eq.~(2.10.18) in \cite{grhyp}.

\begin{Theorem}\label{cnint32}
Let $c=abdefq^{r-1}$. Then there holds
\begin{multline}\label{cnint32gl}
\int_{ax_1}^{b}\dots\int_{ax_r}^{b}
\prod_{1\le i<j\le r}(t_i-t_j)
\prod_{i=1}^r\frac{(qt_i/ax_i,qt_i/b,ct_i;q)_{\infty}}
{(dt_i,et_i,ft_i/x_i;q)_{\infty}}\;
d_qt_r\dots d_qt_1\\
=a^{\binom r2}b^r(1-q)^r
\prod_{1\le i<j\le r}(x_i-x_j)\\\times
\prod_{i=1}^r
\frac{(q,ax_i/b,bq/ax_i,cq^{1-i}/d,cq^{1-i}/e,cx_i/f;q)_{\infty}}
{(adx_i,aex_i,af,bdq^{i-1},beq^{i-1},bf/x_i;q)_{\infty}}.
\end{multline}
\end{Theorem}

Similarly, we can specialize Corollary~\ref{cnnt87} to a multivariable
nonterminating $q$-Pfaff--Saalsch\"utz summation by first replacing
$b$ and $e$ by $aq/b$ and $abq^{1-r}/cdf$, respectively, then letting
$a\to 0$, and finally performing the simultaneous substitutions
$b\mapsto e$, $c\mapsto a$, $d\mapsto b$ and $f\mapsto c$.
We obtain the following multivariable extension of
Eq.~(II.24) in \cite{grhyp}.

\begin{Corollary}[An $A_r$ nonterminating ${}_3\phi_2$
summation]\label{cnnt32}
Let $ef=abcq^r$. Then there holds
\begin{multline}\label{cnnt32gl}
\sum_{S\subseteq\{1,2,\dots,r\}}
\left(\frac qe\right)^{\binom r2-\binom{|S|}2}\prod_{i\notin S}
\frac{(ax_i,bx_i,c,q/ex_i,fq/e;q)_{\infty}}
{(aq/e,bq/e,cq/ex_i,ex_i/q,fx_i;q)_{\infty}}\\\times
\sum_{k_1,\dots,k_r=0}^{\infty}
\underset{i,j\in S}{\prod_{1\le i<j\le r}}
\frac{(x_iq^{k_i}-x_jq^{k_j})}{(x_i-x_j)}
\underset{i,j\notin S}{\prod_{1\le i<j\le r}}
\frac{(q^{k_i}-q^{k_j})}{(x_i-x_j)}
\prod_{i\in S, j\notin S}
\frac{(ex_iq^{1+k_i}-q^{k_j})}{(x_i-x_j)}\\\times
\prod_{i\in S}\frac{(ax_i,bx_i,c;q)_{k_i}}
{(q,ex_i,fx_i;q)_{k_i}}
\prod_{i\notin S}\frac{(aq/e,bq/e,cq/ex_i;q)_{k_i}}
{(q,q^2/ex_i,fq/e;q)_{k_i}}\;\cdot q^{\sum_{i=1}^r k_i}\\
=\prod_{i=1}^r\frac{(q/ex_i,fq^{1-i}/a,fq^{1-i}/b,fx_i/c;q)_{\infty}}
{(aq^i/e,bq^i/e,cq/ex_i,fx_i;q)_{\infty}}.
\end{multline}
\end{Corollary}

It is again clear that Corollary~\ref{cnnt32} above can be
specialized to an $A_r$ terminating $q$-Pfaff--Saalsch\"utz summation,
Namely, by first replacing $f$ by $abcq^r/e$ and then letting $c\to q^{-N}$,
we obtain
\begin{multline}\label{phi32gl}
\sum_{k_1,\dots,k_r=0}^N
\prod_{1\le i<j\le r}\frac{(x_iq^{k_i}-x_jq^{k_j})}{(x_i-x_j)}
\prod_{i=1}^r\frac {(ax_i,bx_i,cx_i,q^{-N};q)_{k_i}}
{(q,cx_i,abx_iq^{r-N}/c;q)_{k_i}}\,\cdot q^{\sum_{i=1} ^rk_i}\\
=\prod_{i=1}^r\frac{(cq^{1-i}/a,cq^{1-i}/b;q)_N}
{(cx_i,cq^{1-r}/abx_i;q)_N},
\end{multline}
which is Theorem~5.1 in \cite{schlhypdet}.

Finally, we specialize Theorem~\ref{cnint32} further, for possible
future reference. We first replace $f$ by $cq^{1-r}/abde$ and then
let $c\to 0$ and replace $a$, $d$, and $e$ by $-a$, $-c/a$, and $d/b$,
respectively. The result is the following.

\begin{Corollary}\label{andaskm}
There holds
\begin{multline}\label{andaskmgl}
\int_{-ax_1}^{b}\dots\int_{-ax_r}^{b}
\prod_{1\le i<j\le r}(t_i-t_j)
\prod_{i=1}^r\frac{(-qt_i/ax_i,qt_i/b;q)_{\infty}}
{(-ct_i/a,dt_i/b;q)_{\infty}}\;
d_qt_r\dots d_qt_1\\
=(-a)^{\binom r2}b^r(1-q)^r
\prod_{1\le i<j\le r}(x_i-x_j)
\prod_{i=1}^r
\frac{(q,-ax_i/b,-bq/ax_i,cdq^{r-1}x_i;q)_{\infty}}
{(cx_i,-adx_i/b,-bcq^{i-1}/a,dq^{i-1};q)_{\infty}}.
\end{multline}
\end{Corollary}

Corollary~\ref{andaskm} is a multivariable extension of a
$q$-integral derived by Andrews and Askey~\cite{andask}.
Replacing $a$, $b$, $c$, $d$, and $x_i$ by
$c$, $d$, $q^a$, $q^b$, and $q^{x_i}$, for $i=1,\dots,r$,
and letting $q\to 1$, we obtain
\begin{multline}\label{mbetangl}
\int_{-c}^d\dots\int_{-c}^d
\prod_{1\le i<j\le r}(t_i-t_j)
\prod_{i=1}^r(1+t_i/c)^{a-1+x_i}(1-t_i/d)^{b-1}\;
dt_r\dots dt_1\\
=\prod_{1\le i<j\le r}(x_i-x_j)
\prod_{i=1}^r
\frac{\Gamma(a+x_i)\,\Gamma(b+i-1)}
{\Gamma(a+b+r-a+x_i)}\\\times
c^{r(1-a)-\sum x_i}\,d^{r(1-b)}\,
(c+d)^{r(a+b-1)+{\binom r2}+\sum x_i},
\end{multline}
which follows from the multiple beta integral evaluation in
Proposition~\ref{mbetaintgl} by the substitutions
\begin{equation}
u_i\mapsto\frac{c+t_i}{d+t_i},\qquad i=1,\dots,r.
\end{equation}


\begin{thebibliography}{99}

\bibitem{andappl} G.~E.~Andrews, Applications of basic hypergeometric
functions, {\em SIAM Rev.\ }{\bf 16} (1974) 441--484.

\bibitem{qandrews} G.~E.~Andrews, {\em $q$-Series: Their development and
application in analysis, number theory, combinatorics, physics and
computer algebra}, CBMS Regional Conference Lectures Series
{\bf 66} (Amer.\ Math.\ Soc., Providence, RI, 1986).

\bibitem{andask} G.~E.~Andrews and R.~Askey, Another $q$-extension of the
beta function, {\em Proc.\ Amer.\ Math.\ Soc.\ }{\bf 81} (1981) 97--100.

\bibitem{sfaar} G.~E.~Andrews, R.~Askey and R.~Roy, {\em Special functions},
Encyclopedia of Mathematics And Its Applications, Vol.~71, Cambridge
University Press, Cambridge, 1999.

\bibitem{bail10} W.~N.~Bailey, An identity involving Heine's basic
hypergeometric series, {J.\ London Math.\ Soc.\ }{\bf 4} (1929) 254--257.

\bibitem{bailwp} W.~N.~Bailey, Well-poised basic hypergeometric series,
{\em Quart.\ J.\ Math.\ }(Oxford) {\bf 18} (1947) 157--166.

\bibitem{bhatdn} G.~Bhatnagar, $D_n$ basic hypergeometric series,
{\em The Ramanujan J.\ }{\bf 3} (1999) 175--203.

\bibitem{degmilne} S. Degenhardt and S.~C.~Milne, A nonterminating
$q$-Dougall summation theorem for hypergeometric series in $U\!(n)$,
in preparation.

\bibitem{grhyp} G.~Gasper and M.~Rahman, {\em Basic hypergeometric series},
Encyclopedia of Mathematics And Its Applications, Vol.~35, Cambridge
University Press, Cambridge, 1990.

\bibitem{geskratt} I.~M.~Gessel and C.~Krattenthaler, Cylindric
Partitions, {\em Trans.\ Amer.\ Math.\ Soc.\ }{\bf 349} (1997) 429--479.

\bibitem{guskradet} R.~A.~Gustafson and C.~Krattenthaler, Determinant
evaluations and $U\!(n)$ extensions of Heine's $_2\phi_1$-transformations,
in: M.~E.~H.~Ismail, D.~R.~Masson and M.~Rahman (Eds.),
{\em Special Functions, $q$-Series and Related Topics},
Amer.\ Math.\ Soc., Providence, R.~I.,
{\em Fields Institute Communications} {\bf 14} (1997) 83--90.

\bibitem{jackson} F.~H.~Jackson, On $q$-definite integrals, {\em
Quart.\ J.\ Pure Appl.\ Math.\ }{\bf 41} (1910) 193--203.

\bibitem{jacksum} F.~H.~Jackson, Summation of $q$-hypergeometric
series, {\em Messenger of Math.\ }{\bf 57} (1921) 101--112.

\bibitem{krattmaj} C.~Krattenthaler, The major counting of
nonintersecting lattice paths and generating functions for tableaux,
{\em Mem.\ Amer.\ Math.\ Soc.\ }{\bf 115}, no.~552 (1995).

\bibitem{milnetf} S.~C.~Milne, Transformations of $U\!(n+1)$ multiple
basic hypergeometric series, in: A.~N.~Kirillov,
A.~Tsuchiya, and H.~Umemura (Eds.), {\em Physics and Combinatorics:
Proceedings of the Nagoya 1999 International Workshop}
(Nagoya University, Japan, August 23--27, 1999), World Scientific,
Singapore, 2001, pp.~201--243. 

\bibitem{milnnew} S.~C.~Milne and J.~W.~Newcomb, Nonterminating
$q$-Whipple transformations for basic hypergeometric series in $U\!(n)$,
in preparation.


\bibitem{schlhypdet} M.~Schlosser, Summation theorems for multidimensional
basic hypergeometric series by determinant evaluations,
{\em Discrete Math.\ }{\bf 210} (2000) 151--169.

\bibitem{schleltf} M.~Schlosser, Elementary derivations of identities
for bilateral basic hypergeometric series, {\em Selecta Math.\ $($N.S.$)$},
to appear.

\bibitem{selberg} A.~Selberg, Bemerkninger om et multiplet integral,
{\em Norske Mat.\ Tidsskr.\ }{\bf 26} (1944) 71--78.

\bibitem{warnell} S.~O.~Warnaar, Summation and transformation formulas
for elliptic hypergeometric series, {\em Constr.\ Approx.\ }{\bf 18}
(2002) 479--502.

\bibitem{whitwatson} E.~T.~Whittaker and G.~N.~Watson, {\em A course of
modern analysis}, 4th ed., Cambridge University Press, Cambridge, 1962.

\end{thebibliography}
\end{document}